\documentclass[ptm,10pt,reqno]{amsart} 
\usepackage{amsmath,amssymb,rawfonts}
\usepackage{tikz}
\usepackage{pgfplots}
\usepackage{color}
\usepackage{mathptmx}
\usepackage[utf8]{inputenc}
\usepackage[all]{xy}
\usepackage{multicol}
\usepackage{hyperref}
\usetikzlibrary{matrix,arrows}

\newtheorem{theorem}{Theorem}
\newtheorem{lem}[theorem]{Lemma}

\begin{document} 

\title{Model Theory, Arithmetic $\&$ Algebraic Geometry}
\author{Joel Torres Del valle}
\date{April, 2019.}
\address{Universidad de Antioquia, Facutad de Ciencias Exactas y Naturales, Instituto de Matemáticas, Medellín-Colombia.}
\email{joel.torres@udea.edu.co}
\vspace{2cm}

\maketitle

\section{Languages and structures.}

\subsection{What is Model Theory?} \hspace{0.2cm}\\

Model Theory is the study of the relationship between theories and mathematical structures which satisfies them \cite{annalissa-toffalori}. Model Theory introduces Mathematical Logic in the practice of Universal Algebra, so we can think it like 

\begin{center}
Model Theory $=$ Universal Algebra $+$ Mathematical Logic. 
\end{center}

I started this pages with the same question Chang-Keisler started their celebrated book \cite{chang-keisler}: `What is Model Theory?' They also summarize the answer to the equation above. In this pages I summarize the relationship between Model Theory, Arithmetic and Algebraic Geometry. The topics will be the basic ones in the area, so this is just an invitation, in the presentation of topics I mainly follow the philosophy of the references, my basic reference is \cite{maker}. The early pioneers in Model Theory were Leopold Löwenheim in 1915, in his work  `Über Möglichkeiten im Relativkalkül'\footnote{Appeared in {\bf Mathematische Annalen} 76: pp. 447–470.} there he proved that if a first order theory $T$ that has an infinite model then has a countable model. Thoralf Skolem in 1920 proved a second version of Löwenheim result, in his work `Logisch-kombinatorische Untersuchungen über die Erfüllbarkeit oder Beweisbarkeit mathematischer Sätze nebst einem Theoreme über dichte Mengen'\footnote{Appeared in {\bf Videnskapsselskapet Skrifter, I. Matematisk-naturvidenskabelig Klasse}. 4: 1-36.}. 

Kurt Gödel in 1930 proved his Completeness and Compactness Theorems for countable languages in his Doctoral dissertation at the University of Viena titled `Die Vollständigkeit der Axiome des logischen Funktionenkalküls'\footnote{Appeared in {\bf Monatshefte für Mathematik} 37: pp. 349-360.}. Alfred Tarski in 1931 proved a third version of Löwenheim result and Anatoly Maltsev in 1936 extended Gödel's Theorems for non-countable languages and showed the existence of non-standard models of Arithmetic by Compactness, a proof which resulted easier than Skolem's one that used Gödel's Imcompleteness Theorems (Gödel proved them in 1931 in his work `Über formal unentscheidbare Sätze dër Principia Mathematica und verwandter systeme I'\footnote{Appearedin {\bf Monatshefte für Mathematik und Physik} 38:  173-198.}; and Skolem proved the existence of non-standard models in 1934 in a work titled `Über die Nicht-charakterisierbarkeit der Zahlenreihe mittels endlich oder abzählbar unendlich vieler Aussagen mit ausschlie$\beta$lich Zahlenvariablen'\footnote{Appeared in {\bf Fundamenta Mathematicae} 1: 150—161.}), Maltsev work was published under the title `Untersuchungen aus dem Gebiete der mathematischen Logik'\footnote{ Appeared in {\bf Matematicheskii Sbornik} 1: 323–336.}.  

More developments were presented by Leon Henkin, Robinson and Tarski in 1940's and the earlies 1950's, see \cite{chang-keisler}, pp. v. Other important developments in 1915-1935 is the method of Quantifiers Elimination, QE to give decision methods for the theories of $(\mathbb{Q}, <)$ by C.H. Langford, and $(\mathbb{Z}, +, <)$ by M. Presburger, and, finally, of the field of complex numbers and of the ordered field of real numbers  Tarski, who proved the completeness of ACF$_{p}$ for $p$ zero or prime, see \cite{casanovas}. Now we introduce languages and semantic of first order logic.  For more see Casanovas's \cite{casanovas}, Villegas-Silva's Book \cite{luis}, Maker's Book \cite{maker} and Chang-Keisler's Book \cite{chang-keisler} also in \cite{page-uno} and \cite{page-dos} are good reads. 

\subsection{Languages} \hspace{0.2cm}\\

\noindent To start we fix a {\bf first order language} $\mathcal{L}$ which contains exactly those symbols that we request in our interest and nothing else. Specifically, this will be given by the union of two big collections: logical constants, containing \emph{connectives}  $\vee$ for \emph{disjunction}, $\wedge$ for \emph{conjunction}, $\to$ for \emph{implication}, $\neg$ for \emph{negation}, $\forall$ for \emph{universal quantifier} and $\exists$ for the \emph{existential quantifier},  a countable collection of \emph{variable symbols} $x_{0}, x_{1}, x_{2},  ...$ and the \emph{equal symbol} $=$. In addition, according to a our interest a set of non-logical constants given by: \emph{function symbols} $f\in \mathcal{F}$ (and an integer  $n_{f}$ for each $f$, called the \emph{range} of $f$), \emph{relation symbols} $R\in\mathcal{R}$ (and an integer  $n_{R}$ for each $R$, called the \emph{arity} of $R$) and  \emph{constant symbols} $c\in \mathcal{C}$, which we assume to have range zero.  

When we say `first order language' is because we are doing first order logic. Basically, this is because our quantification will be on variables. Some good references are Maker's Book \cite{maker} or Tarski's Chapter in \cite{tarski}. A simple example of a language is $\mathcal{L}_{\mathrm{sets}}=\{\in\}$ the language of sets, note that we can define other symbols in Set Theory from $\in$, for example $A\subset B$ iff ($x\in A$ then $x\in B$), or $x\in A\cap B$ iff ($x\in A$ and $x\in B$), etc. In this language we must be careful, we can't say `for all set' in first order logic. On the other hand, see for example that if we consider the language $\mathcal{L}_{a}=\{0, +\}$ the expression $++$ has no meaning but $0+0$ does and also $0+0=0$ but the nature of these two last expressions is different. Now we formalize this distinction.

\subsection{Terms and Formulas of a language} \hspace{0.2cm}\\

\noindent We want to distinguish between chain of symbols defining objects or asserting something  about objects, for example, in $\mathbb{N}$, $0+0$ is an object but $x+x=0$ is an assertion about objects. The {\bf terms} of a language $\mathcal{L}$ will be the smallest set containing all variables and constant symbols in $\mathcal{L}$ and such that: if $t_{1}, ..., t_{n}$ are terms, then $f(t_{1}, ..., t_{n})$ being $f$ a function symbol of range $n$. The formulas of $\mathcal{L}$ are defined by induction: $\varphi$ is an {\it atomic formula} if it is $t_{1}=t_{2}$ or $R(t_{1}, ..., t_{n})$ with $t_{1}, ..., t_{n}$ terms and $R$ a relation symbol. $\varphi$ is a {\bf formula} if it can be obtained by applying a connective or a quantifier on atomic formulas.  

Now, consider the language of rings: $\mathcal{L}_{\mathrm{rings}}=\{0_{R}, 1_{R}, +_{R}, \cdot_{R}, -_{R}\}$ being $0_{R}$ and  $1_{R}$ constants symbols, $+_{R}, -_{R}, \cdot_{R}$ function symbols of range two. Here,  $x+_{R}1_{R}$ is term, while  $x+_{R}1_{R}=0_{R}$ is a formula. Note that the unique variable in  $\psi(x)\equiv\exists x (x+_{R}1_{R}=0_{R})$ in affected by a quantifier. We will omit the subscript $R$ from now on. A formula whose variables are all under the effect of a quantifier is called a {\bf sentence}.  

\subsection{Structures for a language} \hspace{0.2cm}\\

\noindent Let $\mathcal{L}$ be a first order language. A structure for a language $\mathcal{L}$ or simply a $\mathcal{L}$-{\bf structure},  $\mathcal{M}$ with universe  $M$ (a set different of $\emptyset$) is an object of the form  

$$\mathcal{M}=\left(M, \{f^{\mathcal{M}}:M^{n_{f}}\to M\}_{f\in\mathcal{F}}, \{R^{\mathcal{M}}\subset M^{n_{R}}\}_{R\in\mathcal{R}}, \{c^{\mathcal{M}}\in M\}_{c\in\mathcal{C}}\right).$$

That is, a non empty set in which we can give a meaning to each symbol in $\mathcal{L}$. We refer to  $f^{\mathcal{M}}$, $R^{\mathcal{M}}$ and $c^{\mathcal{M}}$ as the \textit{interpretation} of $f, R$ and $c$ respectively.  If $\mathcal{M}$ and $\mathcal{N}$ are structures for $\mathcal{L}$ and $N\subset M$, we say that $\mathcal{N}$ is a {\bf substructure} of $\mathcal{M}$ if operations in $N$ are restrictions of the underlying operations in $M$. For example, an $\mathcal{L}_{\mathrm{rings}}$-structure would be an object $\mathcal{M}=(M, +_{M}, -_{M}, \cdot_{M}, 0_{M}, 1_{M})$. In this case  $+_{M}, -_{M}, \cdot_{M}:M\times M\to M$ are two variable functions  $0_{M}$ and $1_{M}$ are constants. See, for example, that $(\mathbb{Z}, +, -, \cdot, 0, 1)$ and $(\mathbb{Q}, +, -, \cdot, 0, 1)$ are $\mathcal{L}_{\mathrm{rings}}$-structures. 

In this context, given a language, we have some mathematical objects, called structures in which they have a meaning. For example, if we are  given the language of rings, any ring, is an structure for $\mathcal{L} _{\mathrm{rings}}$ but not the converse (in general). 

In mathematical theories given a couple of structures, it is possible in some cases to stablish a comparison criterion by a mapping, for example an homomorphism in Groups Theory, a bijection in Set Theory or a continuous function in Topology (even when this is not a case of the structures we are considering here). Now we generalize this notion to a general structure. 

\subsection{Homomorphisms between structures} \hspace{0.2cm}\\

\noindent Let $\mathcal{M}$ and $\mathcal{N}$ be $\mathcal{L}$-structures, a function $\lambda:\mathcal{M}\to\mathcal{N}$ is called an {\bf homomorphism} provided that

\begin{itemize}
\item[$\bullet$] If $f\in\mathcal{F}$ and $a_{1}, ..., a_{n_{f}}\in M$ then\\ $\lambda(f^{\mathcal{M}}(a_{1}, ..., a_{n_{f}}))=f^{\mathcal{N}}(\lambda(a_{1}), ..., \lambda(a_{n_{f}}))$,
\item[$\bullet$] If $R\in\mathcal{R}$ y $a_{1}, ..., a_{n_{R}}\in M$ then $(a_{1}, ..., a_{n_{R}})\in R^{\mathcal{M}}$ if and only if\\ $(\lambda(a_{1}), ..., \lambda(a_{n_{R}}))\in R^{\mathcal{N}}$ and, 
\item[$\bullet$] If  $c\in\mathcal{C}$ then  $\lambda(c^{\mathcal{M}})=c^{\mathcal{N}}$.
\end{itemize}

If we take a look, the above definition is just what we wanted: a generalization of the notion of homomorphism in Abstract Algebra. For example, consider the groups $\mathbb{Z}$ and $\mathbb{R}$, the mapping $f:\mathbb{Z}\to\mathbb{R}$ defined by $x\mapsto e^{x}$ preserves the interpretation of symbols in the language of groups: $\mathcal{L}_{\mathrm{groups}}=\{*, n\}$ being $*$ a two-adic function symbol and $n$ a constant symbol, in fact, $f(x+y)=e^{x}\cdot e^{y}$ and $f(0)=1$. 

A isomorphism is a bijective homomorphism whose inverse (set-theoretical) is also an homomorphism. Any homomorphism is not necessary a isomorphism, for example, consider for each $a_{1}, ..., a_{n}\in\mathbb{C}$

\[
\begin{array}{ccc}
\Upsilon_{\overline{a}}:(\mathbb{C}[x_{1}, ..., x_{n}], +, -, \cdot, 0, 1)&\longrightarrow&(\mathbb{C}, +, -, \cdot, 0, 1)\\
f(x_{1}, ..., x_{n})&\longmapsto&f(a_{1}, ..., a_{n})
\end{array}
\]

Notice this is an homomorphism of $\mathcal{L}_{\mathrm{rings}}$-structures but not a isomorphism (in fact, $\Upsilon_{a}$ is not injective and hence is not invertible /set-theoretically/). 

\subsection{Truth in a structure}\hspace{0.02cm}\\

\noindent Now, it is a good moment to define what it means that `a sentence $\phi$ is true in a structure $\mathcal{M}$'. Write $t(\overline{a})$ to indicate we are replacing variables in  the term $t$ for $\overline{a}\in M$, and $\phi(\overline{a})$ to indicate we replacing free variables in $\phi$ for $\overline{a}$. We define $\mathcal{M}\models\phi$ by 

\begin{itemize}
\item[$\bullet$] If  $\phi$ is $t_{1}=t_{2}$, then $\mathcal{M}\models\phi(\overline{a})$ if $t_{1}^{\mathcal{M}}(\overline{a})=t_{2}^{\mathcal{M}}(\overline{a})$,
 
\item[$\bullet$] If $\phi$ is $r(t_{1}, ..., t_{n_{r}})$, then $\mathcal{M}\models\phi(\overline{a})$ if  $(t_{1}^{\mathcal{M}}(\overline{a}), ..., t_{n_{r}}^{\mathcal{M}}(\overline{a}))\in R{^\mathcal{M}}$,
 
\item[$\bullet$] If $\phi$ is $\neg\psi$, then $\mathcal{M}\models\phi(\overline{a})$ if $\mathcal{M}\nvDash\psi(\overline{a})$,
 
\item[$\bullet$] If $\phi$ is $(\psi\wedge\theta)$, then $\mathcal{M}\models\phi(\overline{a})$ if $\mathcal{M}\models\psi(\overline{a})$ and $\mathcal{M}\models \theta(\overline{a})$,
 
\item[$\bullet$] If $\phi$ is $\exists v_{j}\psi(\overline{v}, v_{j})$, then  $\mathcal{M}\models\phi(\overline{a})$ if there exits  $b\in M$ such that $\mathcal{M}\models\psi(\overline{a}, b)$,
 
\item[$\bullet$] If $\phi$ is $(\psi\vee\theta)$, then $\mathcal{M}\models\phi(\overline{a})$ if $\mathcal{M}\models\psi(\overline{a})$ or $\mathcal{M}\models\theta(\overline{a})$,
 
\item[$\bullet$] If $\phi$ is $\forall v_{j}\psi(\overline{v}, v_{j})$,  then $\mathcal{M}\models\phi(\overline{a})$ if for all $b\in M$, $\mathcal{M}\models\psi(\overline{a}, b)$. 
\end{itemize}

If $\mathcal{M}\models \phi$ we say that  $\phi$ is \emph{valid} in $\mathcal{M}$. As an example, we see that $\mathbb{R}\models \forall x\hspace{0.1cm}x^{2}\geq0$, $\mathbb{R}
\nvDash\exists x\hspace{0.1cm} x^{2}+1=0$  but $\mathbb{C}\models\exists x\hspace{0.1cm} x^{2}+1=0$ because  $i\in \mathbb{C}$ holds  $i^{2}=-1$, another example is $\mathbb{N}\nvDash\forall x\hspace{0.2cm} c>x$. Latter we will prove that there exits models of arithmetic which satisfies the sentence above, such models are called {\bf non-standard}.   

The fact that we are limited to quantification over elements of the structure is what makes
it `first-order' logic \cite{maker}. 

\subsection{Satisfaction and homomorphism}\hspace{0.02cm}\\

\noindent Let $t(v_{1}, ..., v_{n})$ be a $\mathcal{L}$-term with variables from $v_{1}, ..., v_{n}$, let  $\mathcal{N}$ be a structure and $\mathcal{M}$ a substructure of $\mathcal{N}$ and $b_{1}, ..., b_{n}\in M$. Then $t^{\mathcal{N}}(b_{1}, ..., b_{n})=t^{\mathcal{M}}(b_{1}, ..., b_{n})$, that is, terms are preserved under substructures.  This has an important consequence: if $a_{1}, ..., a_{n}\in M$ with $\phi(v_{1}, ..., v_{n})$ is a quantifier free formula whose free variables are from $v_{1}, ..., v_{n}$. Then 

\begin{center}
$\mathcal{M}\models\phi(a_{1}, ..., a_{n})$ if and only if $\mathcal{N}\models \phi(a_{1}, ..., a_{n})$. 
\end{center}

The proof is easy by induction on formulas, see \cite{maker}. Let  $\mathcal{M}$ and $\mathcal{N}$ be structures for $\mathcal{L}$, we say that these are {\bf elementary equivalent}, and denote it for $\mathcal{M}\equiv\mathcal{N}$ if  $\mathcal{M}\models \phi\Leftrightarrow\mathcal{N}\models \phi$ for any sentence $\phi$. We define  $\mathrm{Th}(\mathcal{M})$ the {\bf complete theory} of  $\mathcal{M}$ as the set of $\mathcal{L}$-sentences $\phi$ such that $\mathcal{M}\models \phi$. 

It holds that $\mathcal{M}\equiv\mathcal{N}$ if and only if $\mathrm{Th}(\mathcal{M})=\mathrm{Th}(\mathcal{N})$.  In fact if  $\mathcal{M}\equiv\mathcal{N}$ then $\phi\in\mathrm{Th}(\mathcal{M})$ implies $\mathcal{M}\models\phi$, as  $\mathcal{M}\equiv\mathcal{N}$ then $\mathcal{N}\models\phi$ so $\phi\in\mathrm{Th}(\mathcal{N})$ for which $\mathrm{Th}(\mathcal{M})\subset\mathrm{Th}(\mathcal{N})$, in the same way if   $\phi\in\mathrm{Th}(\mathcal{N})$ we have $\mathcal{N}\models\phi$ and so  $\mathcal{M}\models\phi$ then  $\phi\in\mathrm{Th}(\mathcal{M})$, hence  $\mathrm{Th}(\mathcal{M})\supseteq\mathrm{Th}(\mathcal{N})$, i.e., $\mathrm{Th}(\mathcal{M})=\mathrm{Th}(\mathcal{N})$. The other sense is obvious.

\begin{lem} \rm Let  $\mathcal{M}$ and  $\mathcal{N}$ be structures. Suppose that  $j:\mathcal{M}\to \mathcal{N}$ is a isomorphism, then if  $t(v_{1}, ..., v_{n})$ is a term with free variables from $v_{1}, ..., v_{n}$ and $a_{1}, ..., a_{n}\in M$ it holds that $j(t^{\mathcal{M}}(a_{1}, ..., a_{n}))=t^{\mathcal{N}}(j(a_{1}), ..., j(a_{n}))$.\end{lem}

The above lemma states just what we expect: terms are preserved under isomorphism. The proof is a simple application of induction on the complexity of terms. The next theorem follows.   

\begin{theorem} \rm Let  $\mathcal{M}$ and $\mathcal{N}$  be structures in a first order language $\mathcal{L}$. Suppose that $j:\mathcal{M}\to \mathcal{N}$ is an isomorphism, then $\mathcal{M}\equiv\mathcal{N}$.\end{theorem}

The proof of the lemma and the theorem above can be found in \cite{maker} and almost in any text about Model Theory. The theorem says, basically that $\mathrm{Th}(\mathcal{M})$ is unique up to isomorphism, i.e., to isomorphic models satisfies exactly the same sentences. It can also be proved that formulas without quantifiers are preserved under sub-structures, that is, if $\mathcal{M}$ and $\mathcal{N}$ are structures such that $\mathcal{M}$ is a sub-structure of $\mathcal{N}$ and $\phi$ is a quantifier free formula of $\mathcal{L}$ then $\mathcal{M}\models\phi$ if and only if $\mathcal{N}\models \phi$, see \cite{maker}.  This remark will be useful below because we will see that theories with QE are model-complete.

\section{ACF and Peano Arithmetic}

\subsection{Fields and Algebraically closed Fields}  \hspace{0.2cm}\\

\noindent A $\mathcal{L}$-{\bf theory} or simply a theory is a set of sentences.  A {\bf model} for a $\mathcal{L}$-theory  $T$ (written $\mathcal{M}\models T$)  is a structure for $\mathcal{L}$ in which each  sentence of $T$ is true, i.e., $\mathcal{M}\models T$ if and only if $\mathcal{M}\models\phi$ for all $\phi\in T$.  

Let FIELDS be is the set of fields axioms

\begin{itemize}
\item[$\bullet$] $\forall x\forall y(x+y=y+x)$ (commutativity of $+$),
\item[$\bullet$] $\forall x\forall y\forall z(x+(y+z)=(x+y)+z)$ (associativity),
\item[$\bullet$] $\forall x(x+0=x)$ (neutral element for addition),
\item[$\bullet$] $\forall x\exists y(x+y=0)$ (inverse element for addition),
\item[$\bullet$] $\forall x\forall y (x\cdot y=y\cdot x)$ (commutativity of product),
\item[$\bullet$] $\forall x\forall y\forall z (x\cdot(y\cdot z)=(x\cdot y)\cdot z)$ (associativity of product),
\item[$\bullet$] $\forall x(x\cdot 1=x)$ (neutral element for product),
\item[$\bullet$] $\forall x(x\neq0\to\exists y(x\cdot y=1))$ (inverse for product),
\item[$\bullet$] $\forall x\forall y\forall z(x\cdot(y+z)=x\cdot y+x\cdot z)$ (distributivity).
\end{itemize}

ACF  algebraically closed fields theory, given by 

\begin{center}
ACF$=$FIELDS+$\forall a_{0}\cdots\forall a_{n}\exists x\left( x^{n}+\displaystyle\sum_{i=0}^{n-1}a_{i}x^{i}=0\right).$
\end{center}

For each integer $n$ and  ACF$_{p}$, the theory of algebraically closed fields of prime characteristic  $p$, given by  ACF$_{p}+\phi_{p}$ for $\phi_{p}$ the formula  

\begin{center}
$\forall x\underbrace{x+\cdots +x}_{p-\mathrm{times}}=0$ and ACF$_{0}=$ACF$+\neg\phi_{n}$ for each integer $n$. 
\end{center}

Any field is a model of FIELDS, any algebraically closed field is a model for ACF. Notice that FIELDS$\subset $ACF, then, any model of ACF will be a model for FIELDS but not the converse. In fact, $(\mathbb{R}, +, -, \cdot, 0, 1)\models$ FIELDS but $(\mathbb{R}, +, -, \cdot, 0, 1)\not\models$ ACF because $x^{2}+1=0$ has not real roots, so $\mathbb{R}\not\models\exists x(x^{2}+1=0)$.

Let  $T$ be a theory, $\mathcal{M}$ and $\mathcal{N}$ model of $T$ such that $N\subset M$. $\mathcal{N}$ is an {\bf elementary sub-structure} of $\mathcal{M}$ if for any formula $\phi(x_{i_{1}}, ..., x_{i_{m}})$ and tuples $a_{i_{1}}, ..., a_{i_{m}}\in N^{m}$, $\mathcal{N}\models\phi(a_{i_{1}}, ..., a_{i_{m}})$ if and only if  $\mathcal{M}\models\phi(a_{i_{1}}, ..., a_{i_{m}})$. In such a case, we write $\mathcal{N}\preceq \mathcal{M}$. 

We say that that  $T$ is \textbf{model-complete} provided thet for any $\mathcal{N}\subset \mathcal{M}$ it holds $\mathcal{N}\preceq \mathcal{M}$. Actually model-completeness is a weak form of QE, that is, a theory with QE is model-complete. The Theory ACF is model-complete, see \cite{maker}. Model completeness was introduced by Robinson, and he showed it could be used to prove Hilbert's Nullstellensatz. 

\subsubsection{Hilbert's Nullstellensatz.}\hspace{0.2cm}\\

Let  $\mathbb{K}\models\mathrm{FIELD}$. Define $\mathbb{K}^{n}=\underbrace{\mathbb{K}\times \cdots\times\mathbb{K}}_{n-\mathrm{times}}$, we call it the  $n$-\textbf{affine space} on $\mathbb{K}$, its elements are called  \textbf{points}. If $f\in\mathbb{K}[x_{1}, ..., x_{n}]$, a point $p=(a_{1}, ..., a_{n})$ in $\mathbb{K}^{n}$ is said a {\bf zero} of $f$ if $f(p)=0$. If $f$ is not identically zero, the set of zeros of  $f$ is called the {\bf hyperfice} defined by  $f$ and is denoted by $V(f)$. 

More generally  $S\subset\mathbb{K}[x_{1}, ..., x_{n}]$, let $V(S):=\{p\in \mathbb{K}^{n}:f(p)=0, p\in S\}$.  We say that  $X\subset \mathbb{K}^{n}$ is a {\bf affine algebraic set}.  Let  $X\subset \mathbb{K}^{n}$, consider those polynomials which are zero on $X$, they form an ideal of  $\mathbb{K}[x_{1}, ..., x_{n}]$ called \textbf{ideal} of $X$, and is denoted by $I(X)$, $I(X)=\{f\in \mathbb{K}[x_{1}, ..., x_{n}]:f(a)=0, \forall a\in X\}$.

\begin{theorem}[Hilbert's Weak Nullstellensatz]\rm
Let  $F$ be an algebraically closed field and $I$ a proper ideal of $F[x_{1}, ..., x_{n}]$. Then $V(I)\neq\emptyset$; i.e., there exits $\overline{a}\in F^{n}$ such that $f(\overline{a})=0$ for all $f\in I$. 
\end{theorem}

This proof of Weak Nullstellensatz will follows from model-completeness of ACF. We will no give a proof now for model-completeness of ACF, it will be in Section 4. From Commutative Algebra we know that $F[x_{1}, ..., x_{n}]$ is Noetherian, then any ideal is finitely generated. Let $f_{1}, ..., f_{k}$ a generator set for $I$ and $J$ a maximal ideal contained in $I$. Define $K=F[x_{1}, ..., x_{n}]/J$. As  $J$ is maximal, $K$ is a field. Even more $F\subset K$.

Let $L$ be the algebraic clausure of  $K$. We have $F\subset K\subset L$, with $F\models$ACF and $L\models$ACF. Note that $\overline{a}:=(x_{1}+M, ..., x_{n}+M)$ is an element in $K^{n}\subset L^{n}$ such that $f_{1}(\overline{a})=\cdots=f_{k}(\overline{a})=0$. Let $b_{1}, ..., b_{M}$ the coefficients in $F$ defining $f_{1}, ..., f_{k}$. Then  

\[
L\models\phi(b_{1}, ..., b_{M}):=\exists\overline{a}\bigwedge_{i=1}^{k}f_{i}(\overline{a})=0
\]

For model-completeness 

\[
F\models\phi(b_{1}, ..., b_{M}):=\exists\overline{a}\bigwedge_{i=1}^{k}f_{i}(\overline{a})=0
\]

Hence there are $\overline{a}\in F^{n}$ such that $f_{1}(\overline{a})=\cdots=f_{k}(\overline{a})=0$. As $f_{1}, ..., f_{k}$ generate $I$, $f(\overline{a})=0$ for any $f\in I$, then  $V(I)\neq\emptyset$. 

\subsection{Peano Arithmetic.} \hspace{0.2cm}\\

Consider the {\bf language of arithmetic} $\mathcal{L}_{A}=\{0, 1, +, -, \cdot, \leq\}$. The Peano Arithmetic is the theory with axioms: associative laws for  $+$ and $\cdot$, its neutral elements $0$ and $1$, respectively, distributive laws of $\cdot$ respect to $+$ and discrete  linear order axioms for $\leq$; 1 is the successor of zero, $\forall x\forall y\forall z (x<y\to x+z<y+z)$, and the induction schema: for any formula $\varphi(x, \overline{w})$, we assume the axiom 

\[
\forall\overline{w}[\varphi(0, \overline{w})\wedge\forall x(\varphi(x, \overline{w})
\to\varphi(x+1, \overline{w}))\to\forall x\varphi(x, \overline{w})].
\]

\noindent $\mathbb{N}\models$PA. Clearly PA$\subset\mathrm{Th}(\mathbb{N})$. Then any model of $\mathrm{Th}(\mathbb{N})$ will be a model of PA,  but not the converse (in general).  We most notice that induction is not an axiom: it is a process which given a formula determines a new axiom when we replace $\varphi$ for a formula in $\mathcal{L}_{A}$, see \cite{maker, Kaye, annalissa-toffalori}.

\subsubsection{Prime models.} \hspace{0.2cm}\\

Let $T$ be a theory with infinite models in a countable language $\mathcal{L}$. We say that $\mathcal{M}\models T$ is a {\bf prime model} provided that given $\mathcal{N}\models T$ there is an elemental embedding from $\mathcal{M}$ in $\mathcal{N}$ (that is, a homomorphism $f:\mathcal{M}\to \mathcal{N}$ such that $f(\mathcal{M})\preceq\mathcal{N}$). Let  $T=$ACF$_{0}$. If $K\models$ACF$_{0}$ and $F=\mathbb{Q}^{\mathrm{alg}}$ then there is an embedding from $F$ in $K$. As ACF$_{0}$ is model-complete this is elemental. So, $F$ is a prime model of ACF$_{0}$. 

Let  $\mathcal{M}\models\mathrm{Th}(\mathbb{N})$, then we can see $\mathcal{N}$ as an initial segment of $\mathcal{M}$, this embedding is elemental, i.e., $\mathbb{N}$ is a prime model of $\mathrm{Th}(\mathbb{N})$, see \cite{maker} for details.

\section{Building Models, Completeness, Compactness and Löwenheim-Skolem Theorems}

\noindent Let $T$ be a first order theory and $\phi$ a sentence, a proof of $\phi$ from $T$ is  a sequence of sentences $\phi_{1}, ...., \phi_{n}$ such that 

\begin{itemize}
\item[$\bullet$] $\phi_{n}=\phi$ and, 
\item[$\bullet$] each $\phi_{i}$ is an axiom or is obtained by applying an inference rule on the others $\phi_{j}$, $1\leq j\leq i$. 
\end{itemize}
 
\noindent A sentence $\phi$ is a \textbf{logical consequence} of a theory $T$ if $\mathcal{M}\models \phi$ for any model $\mathcal{M}$ of $T$, in such a case we write $T\models\phi$, of  $\phi$ in $T$ we write $T\vdash \phi$. The next important result is debt to  Gödel for countable languages in 1930.

\begin{theorem}[Completeness, Gödel-1930.]\rm Let $T$ be a first order theory, then $T\models\phi$ if and only if  $T\vdash\phi$ for any  sentence $\phi$ in the language of  $T$.
\end{theorem} 

This implies that a first order theory has a model if and only if it has not contradictions, that is, there not exits a sentence $\phi$ such that $T\vdash\phi$ and $T\vdash\neg\phi$, that is, a first order theory is satisfacible if and only if is consistent. The following method for building models was introduced by Henkin in his work 

\begin{quote}
`The completeness of the first-order functional calculus'.  Appeared in {\bf Journal of Symbolic Logic}, 14 (1949) pp. 159–166.
\end{quote} 

\noindent To establish the strong completeness of first-order logic. Whilst this method originally involved the deductive apparatus of first-order logic, it can be modified so as to employ only model-theoretic ideas \cite{henkin}.

\subsection{The Compactness Theorem.}\hspace{0.2cm}\\

\noindent A theory  $T$ has {\bf witness property} if for all formula $\phi(v)$ with free variable $v$ there exits a constant $c$ such that $T\vdash\exists v\phi(v)\to\phi(c)$. A theory is {\bf maximal} if for all sentence $\phi$ one $\phi\in T$ or $\neg\phi\in T$. It is a theorem that if $T$ is a finitely satisfacible (i.e., such each finite subset of $T$ is satisfacible) maximal theory  and $T_{i}\subset T$ is finite and $T_{i}\models\psi$ then $\psi\in T$. The proof is simple by contradiction. If we suppose $\psi\not\in T$, then $\neg\psi\in T$ because $T$ is maximal \cite{maker}.

The above means that $T_{j}\subset T$ finite is such that $T_{j}\models\neg\psi$ but there exits $T_{i}\subset T$ finite such that $T_{i}\models\psi$, so we can find a finite subset of $T$ which has no models and that is impossible because $T$ is finitely satisfacible. A finitely satisfacible theory $T$ with the witness property always has a model, see \cite{maker} for a proof. 

Given a finitely satisfacible theory  $T$, it  is possible to find  a language $\mathcal{L}^{*}$ containing $\mathcal{L}$ (the language of $T$) and a $\mathcal{L}^{*}$-theory $T^{*}$ containing $T$ finitely satisfacible such that any $\mathcal{L}^{*}$-theory extending $T^{*}$ which has the witness property. We can take $\mathcal{L}^{*}$ such that $|\mathcal{L}^{*}|=|\mathcal{L}|+\aleph_{0}$. Note that this model will also be a model of $T$.  For $T$ a $\mathcal{L}$-theory finitely satisfacible and $\phi$ a $\mathcal{L}$-sentence, then any $T\cup\{\phi\}$ or $T\cup\{\neg\phi\}$ is finitely satisfacible. Then any $\mathcal{L}$-theory finitely satisfacible $T$ has maximal extension finitely satisfacible. If $T$ is a $\mathcal{L}$-theory finitely satisfacible and $\kappa$ an infinite cardinal holding $\kappa\geq|\mathcal{L}|$, there is $\mathcal{L}^{*}$ containing $\mathcal{L}$ and $T^{*}$ containing $T$ finitely satisfacible with the witness property $|\mathcal{L}^{*}|$ at most $\kappa$. There is a $\mathcal{L}^{*}$-theory $T'$ extending  $T^{*}$ finitely satisfacible and maximal. As $T'$ has the witness property there exits a model $\mathcal{M}$ of $T$ whose cardinal is at most $\kappa$ \cite{maker, chang-keisler}. 
 
\begin{theorem}[Compactness Theorem, Gödel-Maltsev 1930-1936.] \rm If $T$ is a $\mathcal{L}$-theory  finitely satisfacible and $\kappa$ a infinite cardinal with $\kappa\geq|\mathcal{L}|$, then there exists a model of $T$ with cardinal at most $\kappa$.  \end{theorem}
 
The Compactness Theorem was first proved by Gödel in 1930 as a consequence of the Completeness Theorem. This result was shown to be true for countable languages. Maltsev extended this result to non countable languages in 1936.

Now we will see some implications of Compactness Theorem. 

Consider $T$ the theory obtained by joining to FIELDS the following axioms  
 
$\neg(1+1=0)$

$\neg(1+1+1=0)$

\hspace{1cm}\vdots

$\neg(\underbrace{1+1+1+...+1}_{n-\mathrm{times}}=0)$

\hspace{1cm}\vdots

If $\varphi$ is any sentence true in any fields of characteristic zero,  $T\cup\{\neg\varphi\}$ has no models, and then it has a finite subset which has no models. Hence there exits  $n\in\mathbb{N}$ such that $\varphi$ is true in any fields of characteristic bigger than  $n$.   We had proved that a sentence is true in any field of characteristic zero if and only if it is true in any field of characteristic $p$ for $p\geq n$ for a fix $n$ enough large, \cite{keisler}. 

Another application is the existence of non-standard models of $\mathrm{Th}(\mathbb{N})=\{\phi:\mathbb{N}\models\phi\}$,  that is to say,  models not isomorphic to $\mathbb{N}$. 

Consider $T$ the theory obtained by join to $\mathrm{Th}(\mathbb{N})$ the axioms $c>\underbrace{1+\cdots+1}_{n-\mathrm{times}}$ for all $n\in\mathbb{N}$ and being $c$ a new constant added to the language $\mathcal{L}_{A}$. The theory $T$ is finitely satisfacible and then satisfacible. For that we just note that a finite subset $T_{0}$ of $T$ in a theory of the form $\mathrm{Th}(\mathbb{N})\cup T^{k}$ being $T_{k}=\{c>\underbrace{1+\cdots+1}_{n-\mathrm{times}}\}_{n\leq k}$, and it is clear that $\mathbb{N}\models \mathrm{Th}(\mathbb{N})\cup T^{k}$. Then, there exits a model for $T$. That is to say, there exits  $c\in\mathcal{M}$ such that $c>n$ for all $c\in\mathbb{N}$. This models can't be isomorphic to $\mathbb{N}$, they are called {\bf non-standard models of arithmetic}, see Kaye's Book \cite{Kaye}.

\subsection{Löweheim-Skolem-Tarski Theorem.}\hspace{0.2cm}\\

\begin{theorem}[Löwenheim-Skolem-Tarski]\rm  Let $T$ be a first order theory in the language $\mathcal{L}$. If $T$ has at least an infinite model, then it has a model of in infinite cardinal.\end{theorem}

That is to say, first order theories can't control the cardinal of their models. If any two models of $T$ with cardinal $\kappa$ are isomorphic, we will say that $T$ is $\kappa$-{\bf categoric}. A consequence of the theorem above is that non-standard models of PA and $\mathrm{Th}(\mathbb{N})$ can be taken of cardinal as big a we prefer, so there are infinitely many non-isomorphic models of Arithmetic. Certainly, there exits exactly $2^{\aleph_{0}}$ non-isomorphic models of $\mathrm{Th}(\mathbb{N})$, see \cite{Kaye}.

If $T$ is a theory without finite models, and $T$ is $\kappa$-categoric, then $T$ is {\bf complete}, that is, for any sentence of the language of $T$, $T\models\phi$ or $T\models\neg\phi$, this result is known as {\it Vaught's Test}. The proof is again simple by contradiction, if $T$ was not complete then for some $\mathcal{L}$-sentence $\phi$ we must have $T\not\models\phi$ and $T\not\models\neg\phi$, the theories $T\cup\{\phi\}$ and $T\cup\{\neg\phi\}$ are satisfacible.  Because $T$ has not finite models both  $T\cup\{\phi\}$ and $T\cup\{\neg\phi\}$ have infinite models, say, $\mathcal{M}_{\phi}$ and $\mathcal{M}_{\neg\phi}$, respectively and both can be of cardinality $\kappa$. As $T\cup\{\phi\}$ and  $T\cup\{\neg\phi\}$ differs respect to $\phi$ we don't have elementary equivalence (recall the last part of Section 1), so we don't have isomorphism, a contradiction. 

Morley proved in 1965 that if a first order theory $T$ is $\kappa$-categoric for $\kappa$ non-countable, then for any $\lambda$ non-countable $T$ is $\lambda$-categoric. This result had been conjectured by Lo\'s in 1954. 

Morley result in contrasting with Löwenheim-Skolem-Tarski one, can be seen as follows: a first order theory $T$ which is satisfacible has a model in any infinite cardinal. But, if $T$ is $\kappa$-categoric for some infinite cardinal $\kappa$, we can see that up to isomorphism $T$ has exactly one model of cardinal $\kappa$. Morley result then says that a $\kappa$-categoric theory has one model of cardinal a fix infinite cardinal $\lambda$. So, models of $T$ a determined up to isomorphism by its cardinal if $T$ is $\kappa$-categoric. 

ACF$_{p}$ is $\kappa$-categoric and hence complete, for $p$ zero or prime. In fact, if $A$ and $B$ be two models of ACF$_{p}$, then they are infinite and if they have the same cardinality, they most be isomorphic because any algebraically closed field is unique determined by the its prime characteristic and its transcendence degree by the Steinitz's classification Theorem for ACF. Then $T$ is $\kappa$-categoric and hence complete. 

\begin{lem}[Lefschetz's Principle]\rm
Let $\varphi$ be a formula in $\mathcal{L}_{R}$, the following assertions are equivalent

\begin{itemize}
\item[$\bullet$] For enough big primes $p$, there is a model $K\models ACF_{p}\cup\{\varphi\}$;
\item[$\bullet$] ACF$_{0}\models\varphi$; 
\item[$\bullet$] $\mathbb{C}\models\varphi$.
\end{itemize} 
\end{lem}

The proof of this lemma follows easily by the application of the ideas above, for example, the last two item are equivalent for completeness of ACF$_{0}$. The equivalence of the two first was shown before with the Compactness Theorem. 

\subsubsection{Noether-Ostroski Irreducibility Theorem} \hspace{0.2cm}\\

Given $f\in\mathbb{Z}[x_{1}, ..., x_{n}]$ fix and $p$ prime. Let $f_{p}\in\mathbb{F}_{p}[x_{1}, ..., x_{n}]$ obtained from  $f$ by reducing any coefficient module $p$. 

\begin{theorem}[Noether-Ostroski Irreducibility Theorem]\rm Let $f\in\mathbb{Z}[x_{1}, ..., x_{n}]$ fix and $p$ prime.  $f$ is irreducible in $\mathbb{Q}^{\mathrm{alg}}$ if and only if $f_{p}$ is irreducible in $\mathbb{F}^{\mathrm{alg}}_{p}$ for co-finitely many primes. \end{theorem}
 
We know $\mathbb{Q}^{\mathrm{alg}}\models ACF_{0}$. Let $n=\deg(f)$, given  $k\in\mathbb{N}$, define $f_{\overline{a}^{k}}$ the general polynomial of grade $k$  generated by the coefficient $\overline{a}^{k}:=(a_{1}^{k}, ..., a_{n_{k}}^{k})$. Define $\phi$ by  
  
$$\phi:=\bigwedge_{k+l=n}\forall\overline{a}^{k}\forall\overline{b}^{l}(f_{\overline{a}^{k}}f_{\overline{b}^{k}}\neq f).$$
  
\noindent In this case $\mathbb{Q}^{\mathrm{alg}}\models\phi$. By the previous lemma we have the result.

\section{Quantifier elimination, QE}

Let $\phi$ and $\psi$ two formulas in  $\mathcal{L}$, it could occur that they admit a same meaning in a structure $\mathcal{M}$ or a class of structures. Even more,  that a theory  $T$ satisfies $\phi(\overline{v})\leftrightarrow\psi(\overline{v})$. For example, the formula  $\varphi(v):v\geq0$ has the same meaning of  $\psi(v):\exists w(v=w^{2})$ in $\mathbb{R}$ or, in the integers  $\exists w_{1}\exists w_{2}\exists w_{3}\exists w_{3} (v=w_{1}^{2}+w_{2}^{2}+w_{3}^{2}+w_{4}^{2})$ \cite{annalissa-toffalori, keisler}.

In particular, we are interested in theories for which given a formula  $\phi(\overline{v})$ it is possible to find a formula in which there are not quantifiers such that it has the same meaning of $\phi(\overline{v})$. A such theory is said to has {\bf quantifier elimination}.  

In general, the method for proving quantifiers is the following: first, given a theory $T$ we get a set of formulas called basic formulas. For a Boolean combination we understand a formula obtained by using connectives in basic formulas, we show that any formula can be written as a boolean combination of basic formulas. 

A useful theorem for proving quantifier elimination is the following. 

\begin{theorem}[Quantifier Elimination, QE]\rm $T$ has quantifier elimination if $\mathcal{M}, \mathcal{N}\models T$, $\mathcal{A}$ is a $\mathcal{L}$-structure, $A\subset M$ and $A\subset N$, then $\mathcal{M}\models\phi(a_{1}, .., a_{n})$ if and only if $\mathcal{N}\models\phi(a_{1}, ..., a_{n})$ for all $a_{1}, ..., a_{n}\in A$.
\end{theorem} 

\noindent Notice that formulas without quantifiers are preserved under sub-structures, then theories with QE are model-complete. A set $X$ is said to be {\bf definable} if there exists a formula $\phi(v_{1}, ..., v_{n}, w_{1}, ..., w_{m})$ and $(b_{1}, ..., b_{m})\in M^{m}$ such that  $X=\{(a_{1}, ..., a_{n})\in M^{n}:\mathcal{M}\models\phi(a_{1}, ..., a_{n}, b_{1}, ..., b_{m})\}$. We say that  $\phi$ {\bf defines} $X$. Now, let us see some classical examples of definable sets.  

Robinson proved that  

$$\forall y\forall z ([\phi(y, z, 0)\wedge(\forall w(\phi(y, z, w)\to\phi(y, z, w+1)))]\to \phi(t, z, x))$$

\noindent defines integer in $\mathbb{Q}$, being $\phi$ the formula $\exists a\exists b\exists c xyz^{2}+2=a^{2}+xy^{2}-yc^{2}$. For the remarks above, we can also define positive integers in the integers. But we can't define real numbers in complex numbers  \cite{maker, annalissa-toffalori}. 

ACF has QE an also ACF$_{p}$ with $p$ zero or prime. The completeness  of ACF$_{p}$ for $p=0$ or prime was first proved by Alfred Tarski using the method of QE. Actually he gave an explicit algorithm for eliminating quantifiers.  Any theory with QE is model-complete, so it shows model-completeness for ACF and ACF$_{P}$.

Quantifiers Elimination permits a geometrical interpretation. A theory is said to be {\bf strongly minimal} if for any $\mathcal{M}\models T$ and $X\subset M$ is definable, then  $X$ or its complement is finite. ACF is strongly minimal for if $A$ a definable subset of $F\models$ACF, $A$ is a Boolean combination of set of the form  $\{x:f(x)=0\}$ with $f\in F[\overline{X}]$. If $f$ is not identically zero, then $V_{f}(\{f\})$ is finite. Actually, definable sets are exactly those constructible from algebraic geometry.

Let $F\models$ FIELDS and $X\subset F^{n}$, we say that  $X$ is {\bf Zariski closed} if it is a Boolean combination of set $\left\{\overline{x}:\displaystyle\bigwedge_{i=1}^{m}f_{i}(\overline{x})=0\right\}$ with  $f_{1}, ..., f_{m}\in F[\overline{X}]$. By Hilbert's Basis Theorem, the intersection of Zariski closed is again a Zariski closed, for which they form a topology for  $F^{n}$.  A set is called {\bf constructible} if it is a Boolean combination of  Zariski closed, see \cite{maker}.

\noindent Let $\mathcal{L}$ be a language, $\mathcal{M}$ a model for $\mathcal{L}$ and $\mathcal{L}_{A}$ the language resulting of adding a new constant symbol $\underline{a}$ to $\mathcal{L}$ for each $a\in A\neq\emptyset$, let $p$ be a set of  $\mathcal{L}_{A}$-formulas in free variables $v_{1}, ..., v_{n}$. We say that  $p$ is a $n$-{\bf type} if  the theory 

\begin{center}
$p\cup\mathrm{Th}_{A}(\mathcal{M})=\{\phi:\mathcal{M}\models\phi$ with $\phi$ a $\mathcal{L}_{A}$-formula$\}$ is satisfacible.
\end{center}

\noindent If for all  $\mathcal{L}_{A}$-formula $\phi$ with free variables from $v_{1}, ..., v_{n}$ it holds $\phi\in p$ or $\neg\phi\in p$ then $p$ is said to be \textbf{complete}. We denote $S^{\mathcal{M}}_{n}(A)$ the set of complete $n$-types. Let $\overline{\phi}:=\{p\in S_{n}^{\mathcal{M}}(A):\phi\in p\}$. 

It is easy to prove that $\overline{\phi\vee\psi}=\overline{\phi}\cup\overline{\psi}$ and $\overline{\phi\wedge\psi}=\overline{\phi}\cap\overline{\psi}$. We consider {\bf Stone's topology} on $S_{n}^{\mathcal{M}}(B)$ by taking basic open set those of the form $\overline{\phi}$. The topological space of types was introduced by Morley in \cite{morley} solving a conjecture of Lo\'s. Note that $\overline{\phi}=S^{\mathcal{M}}_{n}(A)-\overline{\neg\phi}$ so, any $\overline{\phi}$'s are clopen. The mapping  

\[
\begin{array}{ccc}
\Psi:S_{n}^{K}(k)&\longrightarrow&\mathrm{Spec}(K[X_{1}, ..., X_{n}])\\
p&\longmapsto&I_{p}:=\{f\in F[X_{1}, ..., X_{n}]:f(\overline{x})=0\in p\}
\end{array}
\]

\noindent is a continuous bijection: suppose that  $f_{1}, ..., f_{m}\in k[X_{1}, ..., X_{n}]$, then 

\[
\Psi^{-1}(\{P\in \mathrm{Spec}(k[X_{1}, ..., X_{n}]):f_{1}, ...., f_{m}\in P\})=\left\{p\in S_{n}^{K}(k):\bigwedge_{i=1}^{m} f_{i}(\overline{v})=0\in p\right\}
\]

\noindent a clopen set. So,  $p\mapsto I_{p}$ is continuous. For biyectivity see \cite{maker}.  The mapping  $p\to I_{p}$ is not a homeomorphism. In particular, if $K$ is an algebraically closed field and  $k$ a sub-field of $K$,  $f\in k[X_{1}, ..., X_{n}]$ is not constant, $\{p\in S_{n}^{K}(k):f(\overline{v})=0\}$ is clopen in  $S_{n}^{K}(k)$, while the image in  $\mathrm{Spec}(A)$ is closed but not open. 

The Stone space $S_{n}^{\mathcal{M}}(A)$ is compact \cite{maker}. Zariski topology on $\mathrm{Spec}(k[X_{1}, ..., X_{n}])$ is compact, because $\mathrm{Spec}(k[X_{1}, ..., X_{n}])$  is the image continuous of a compact, that is to say, $S_{n}^{K}(k)$, see \cite{maker}.  By Hilbert Basis Theorem, any ideal in $k[X_{1}, ..., X_{n}]$ is finitely generated. So, there just exits $|k|+\aleph_{0}$ prime ideals on $K\models$ACF and $k$ a sub-field of $K$. Then, there exists exactly $|k|+\aleph_{0}$ $n$-complete types on $k$. 

Given a theory $T$ the best possible situation is $\kappa$-categoricity, but this is not always the case. So, we need to explore similar (but weaker) conditions which permits us to classify models of a first order theory. This notion will be {\bf stability}, remember that Löwemhein-Skolem-Tarski Theorem proves that a first order theory has a model in any infinite cardinal. A theory $T$ is $\omega$-\textbf{stable} if for any model $F$ of $T$, $|S_{n}(T)|=|F|$, this means, grosso modo, the theory $T$ has not many types. 

Remember that  $|\mathrm{Spec}(F[X_{1}, ..., X_{n}])|=|F|$, $F\models$ACF, $|S_{n}(F)|=|F|$, because the mapping $S_{n}(F)\to \mathrm{Spec}(F[X_{1}, ..., X_{n}])$ applying $p\mapsto I_{p}$ is bijective, then ACF$_{p}$ is $\omega$-stable for $p$ zero or prime. 

In 1960, R. Vaught conjectured a complete first-order theory has either countable many or countable models. Sharon Shelah established the result for $\omega$-stable theories. Buechler and Newelski extended the result to certain superstable theories (i.e., stable /that is, has exactly $\lambda$ complete types, for $\lambda$ a cardinal/ in all cardinalities beyond the continuum). 

Stability theory interacts with classical Algebra in several ways. In \cite{nee} it is shown that every stable semi-simple ring is a matrix ring over an algebraically closed field; in \cite{che} it is shown that every superstable division ring is an algebraically closed field, see \cite{page-dos}. 

This last discussion is based in \cite{page-dos} and \cite{page-uno}. Other important growing area is {\bf Zariski Geometries} of Zilber, B; and  Hrushovski, E., the main reference is

\begin{quote}
`Zariski Geometries: geometry from the logicians point of view'. London Mathematical Society, \textbf{Lectures Series 360}. Cambridge University press. By Boris Zilber. 
\end{quote}

And the original work {\it Zariski Geometries}, \textbf{Bulletin of the American Mathematical Society}. 	Vol. 28, No. 2, 1993. pg. 315-323 by Zilber, B; and  Hrushovski, E. Their main theorem is: Let $X$ be a very ample Zariski geometry. Then there exists a smooth curve $C$ over an algebraically closed field $F$, such that $X$ and $C$ are isomorphic as Zariski geometries. $F$ and $C$ are unique determined, up to a fields isomorphism and a isomorphism of curves over $F$, see \cite{zilber} for details.

\end{document}